\documentclass[12pt]{amsart}
\usepackage{amsmath,amsfonts,amssymb}
\numberwithin{equation}{section}
\newtheorem{theorem}{Theorem}[section]
\newtheorem{corollary}[theorem]{Corollary}

\newtheorem{proposition}[theorem]{Proposition}
\newtheorem{remark}[theorem]{Remark} 
\def\a{\mathfrak{a}}
\def\h{\mathfrak{h}}

\def\R{\mathbb{R}}
\def\C{\mathbb{C}}

\title{
A formula for the hypergeometric function of type $BC_n$
}

\author{Nobukazu Shimeno}
\subjclass[2000]{ 
Primary 33C67; Secondary 43A90.}
\address{Department of Applied Mathematics, 
Okayama University of Science, 
Okayama 700-0005, Japan}
\email{shimeno@xmath.ous.ac.jp}
\date{}
\begin{document} 
\maketitle

\begin{abstract}
Formulae of Berezin and Karpelevi\v c for the radial parts of 
invariant differential operators and the spherical function 
on a complex Grassmann manifold are generalized to the 
hypergeometric functions associated with root system of 
type $BC_n$ under condition that the multiplicity of the 
middle roots is zero or one. 
\end{abstract}

\section*{Introduction}
Berezin and Karpelevi\v c \cite{BK} gave an explicit expression for 
radial parts of invariant differential operators and 
the spherical functions on $SU(p,q)/S(U(p)\times U(q))$ without proof and 
Hoogenboom~\cite{Hoogenboom} gave proof of these results. 
Explicit expressions of the Laplace-Beltrami operator and  higher order 
invariant differential operators allows us to construct 
eigenfunctions by the method of separation of variables 
and the spherical function has an expression using determinant of a matrix 
whose entries are the Gauss hypergeometric functions. 

Heckman and Opdam developed theory of the hypergeometric 
function associated with a root system, 
which is a generalization 
of theory of the spherical function on a symmetric space (\cite{Heckman}). 
Namely, the radial part of the Laplace-Beltrami operator of a 
Riemannian symmetric space of the non-compact type consists 
of data such as the restricted root system, multiplicities of 
roots. Heckman and Opdam allowed multiplicities of roots 
arbitrary complex numbers (that coincide on every Weyl group orbit) 
and constructed commuting family of differential operators and 
eigenfunctions. For rank one (one variable) case, 
their hypergeometric function is the Jacobi function (\cite{Koornwinder}), 
which is essentially the same as the Gauss hypergeometric function. 

In this paper, the author proves that 
the results of Berezin and Karpelevi\v c \cite{BK} 
are valid for the hypergeometric function associated with root 
system of type $BC_n$ under the condition that the multiplicity 
of the middle roots is $1$. Though it is an easy generalization of \cite{BK}, 
our results cover integral middle multiplicities in conjunction with the 
hypergeometric shift operator, which include many cases of symmetric spaces. 

\section{Hypergeometric function associated with a root system}
\subsection{Notation}
In this section, we review on the hypergeometric function 
associated with a root system. 
See \cite{Heckman} for details. 

Let $E$ be an $n$-dimensional 
Euclidean space with inner product $(\cdot,\cdot)$. 
For $\alpha\in E$ with $\alpha\not=0$ write
\begin{equation}
\alpha^\vee=\frac{2\alpha}{(\alpha,\alpha)}.
\end{equation}
Let $R\subset E$ be a root system of rank $n$ and $W$ its Weyl 
group. Let $R_+\subset R$ be a fixed set of positive roots and 
$E_+\subset E$ be the corresponding positive Weyl chamber. 
Let 
\begin{equation}
P=\{\lambda\in E\,:\,(\lambda,\alpha^\vee)\in\mathbb{Z}\,\,
\forall\,\alpha\in R\}.
\end{equation}

Let $k_\alpha\,(\alpha\in R)$ be complex numbers such that 
$k_{w\alpha}=k_\alpha$ for all $w\in W$. We call $k=(k_\alpha)_
{\alpha\in R}$ a multiplicity function on $R$. 
Let $K$ denote the set of multiplicity function on $R$. 
We put 
\begin{align}
& \rho(k)=\frac12\sum_{\alpha\in R_{+}}k_\alpha \alpha, \\
& \delta(k)=\prod_{\alpha\in R_+}(e^{\frac12\alpha}
-e^{-\frac12\alpha})^{2k_\alpha}.
\end{align}

\subsection{Commuting family of differential operators}
Let $\xi_1,\dots,\xi_n$ be an orthonormal basis of $E$ and 
consider the differential operator
\begin{equation}
\label{eqn:laplace1}
L(k)=\sum_{j=1}^n \partial_{\xi_j}^2+\sum_{\alpha\in R_+}
k_\alpha\frac{1+e^{-\alpha}}{1-e^{-\alpha}}
\partial_\alpha
\end{equation}
on $E$. Here $\partial_\alpha$ denotes the directional 
derivative along $\alpha$ such that $\partial_\alpha (e^\lambda)
=(\alpha,\lambda)e^\lambda$ for $\alpha,\,\lambda\in E$. 
We have
\begin{align}
\label{eqn:laplace2}
\delta(k)^{\frac12} & 
\circ\{L(k)+(\rho(k),\rho(k))\}\circ\delta(k)^{-\frac12} \\
& = \sum_{j=1}^n \partial_{\xi_j}^2+\sum_{\alpha\in R_+}
\frac{k_\alpha(1-k_\alpha-2k_{2\alpha})(\alpha,\alpha)}
{(e^{\frac12\alpha}-e^{-\frac12\alpha})^2}.\notag
\end{align}

Let $\mathcal{R}$ denote the algebra generated by the functions
\begin{equation}
\frac{1}{1-e^{-\alpha}}\quad (\alpha\in R_+)
\end{equation}
viewed as a subalgebra of the quotient field of $\mathbb{R}[P]$. 
Let $S(E)$ denote the symmetric algebra of $E$. 
Let $\mathbb{D}_\mathcal{R}=\mathcal{R}\otimes S(E)$ denote 
the algebra of differential operators on $E$ with coefficient in 
$\mathcal{R}$ and let $\mathbb{D}_\mathcal{R}^W$ be the 
subalgebra of $W$-invariants in $\mathbb{D}_\mathcal{R}$. 
Let $\gamma(k)$ denote the algebra homomorphism 
\begin{equation}
\gamma(k)\,:\,\mathbb{D}_\mathcal{R}\longrightarrow S(E)
\end{equation}
defined by 
\[
\gamma(k)\left(\frac{1}{1-e^{-\alpha}}\right)=1\quad (\alpha\in R_+).
\]
Let 
\begin{equation}
\mathbb{D}(k)=\{D\in \mathbb{D}_\mathcal{R}^W \,:\, 
[L(k),P]=0\}
\end{equation}
denote the commutator of $L(k)$ in $\mathbb{D}_\mathcal{R}^W$ 
and let $S(E)^W$ denote the set of $W$-invariants in $S(E)$. 

\begin{theorem}
\label{thm:ho1}
The map 
\begin{equation}
\gamma(k)\,:\,\mathbb{D}(k)\longrightarrow S(E)^W
\end{equation}
is an algebra isomorphism. In particular, $\mathbb{D}(k)$ is a 
commutative algebra. Moreover, 
if $D\in \mathbb{D}_\mathcal{R}^W$ is a differential operator of 
order $N$, then its principal symbol $\sigma(D)$ has constant coefficients 
and coincides with homogeneous component of $\gamma(k)(D)$ of 
degree $N$. 
\end{theorem}

\subsection{The hypergeometric function}
Let $Q$ be the root lattice
\begin{equation}
Q=\{\textstyle
\sum_{\alpha\in R_+}z_\alpha \alpha\,:\,z_\alpha\in \mathbb{Z}_+\}. 
\end{equation}
Put
\begin{equation}
\mathfrak{h}=E_\mathbb{C}=\mathbb{C}\otimes_\mathbb{R} E,
\quad 
A=\exp\,E,\quad e=\exp 0,\quad A_+=\exp E_+
\end{equation}
For $\mu\in\h^*$ and $a\in A$, we write $a^\mu=\exp(\mu(\log a))$. 

If $\lambda\in\h^*$ satisfies the condition
\begin{equation}
\label{eqn:generic}
-2(\lambda,\mu)+(\mu,\mu)\not=0\text{ for all }\mu\in Q,
\end{equation}
then the equation
\begin{equation}
\label{eqn:hgl}
L(k)u=((\lambda,\lambda)-(\rho(k),\rho(k))u
\end{equation}
has a unique solution on $A_+$ of the form
\begin{equation}
\label{eqn:hcseries}
u(a)=\Phi(\lambda,k;a)=\sum_{\mu\in Q}\Gamma_\mu a^{\lambda-\rho(k)-\mu}
\end{equation}
with $\Gamma_0=1$. The function $\Phi(\lambda,k;a)$ is also a 
solution of the system of differential equations
\begin{equation}
\label{eqn:hgs}
D u=\gamma(k)(D)(\lambda)u,\quad D\in\mathbb{D}(k).
\end{equation}
If 
\[
(\lambda,\alpha^\vee)\not\in\mathbb{Z}\text{ for all }
\alpha\in R,
\]
then $\Phi(w\lambda,k;a)\,\,(w\in W)$ form a basis of the 
solution space of (\ref{eqn:hgs}). 

Define meromorphic functions $\tilde{c}$ and $c$ on $\h\times K$ by
\begin{equation}
\label{eqn:cf1}
\tilde{c}(\lambda,k)=\prod_{\alpha\in R_+}
\frac{\Gamma((\lambda,\alpha^\vee)+\frac12k_{\frac12\alpha})}
{\Gamma((\lambda,\alpha^\vee)+\frac12k_{\frac12\alpha}+k_\alpha)}
\end{equation}
and
\begin{equation}
\label{eqn:cf2}
c(\lambda,k)=\frac{\tilde{c}(\lambda,k)}{\tilde{c}(\rho(k),k)}
\end{equation}
with the convention $k_{\frac12\alpha}\not=0$ if $\frac12\alpha\not\in R$. 
We call the function
\begin{equation}
F(\lambda,k;a)=\sum_{w\in W}c(w\lambda,k)\Phi(w\lambda,k;a)
\end{equation}
the hypergeometric function associated with $R$. 
Let $S\subset K$ denote the set of zeroes of $\tilde{c}(\rho(k),k)$. 

\begin{theorem}
\label{thm:ho2}
Assume that $k\in K\setminus S$. Then the system of 
differential equation {\em (\ref{eqn:hgs})} has a unique solution 
that is regular at $e\in A$, $W$-invariant, and 
\[F(\lambda,k;e)=1.\]
The function $F$ is holomorphic in $\lambda\in \h,\,k\in K\setminus S$, 
and analytic in $a\in A$.  
\end{theorem}

\begin{remark}
\label{remark:sym}
{\em 
Theorem~\ref{thm:ho1} and Theorem~\ref{thm:ho2} 
were proved by Heckman and Opdam 
in a series of papers. See \cite{Heckman} and references therein. 

Let $G/K$ be a Riemannian symmetric space of the non-compact type, 
$\Sigma$ be the restricted root system, and $m_\alpha$ be 
the root multiplicity (dimension of the 
root space) of $\alpha\in \Sigma$. 
Put 
\begin{equation}
R=2\Sigma,\quad k_{2\alpha}=\frac12m_\alpha.
\end{equation}
Then {(\ref{eqn:laplace1})} is the radial part of the Laplace-Beltrami 
operator on $G/K$, $\mathbb{D}(k)$ is the algebra of radial parts 
of invariant differential operators on $G/K$, 
and $F(\lambda,k;a)$ is the radial part of the spherical function 
on $G/K$. 
In this case Theorem~\ref{thm:ho1} and Theorem~\ref{thm:ho2} 
were previously proved by Harish-Chandra. 
See \cite{Helgason} for theory of spherical functions 
on symmetric spaces. }
\end{remark}

\subsection{Rank one case}
For a root system of rank $1$, the hypergeometric function 
is given by the Jacobi function. We review on the Jacobi function. 
See \cite{Koornwinder} for details. 

Assume that $R=\{\pm e_1,\,\pm 2e_1\}$ 
with $(e_1,e_1)=1$ and 
put 
\begin{equation}
\label{eqn:alpha}
k_s=k_{e_1},\quad k_l=k_{2e_1},\quad 
\alpha=k_s+k_l-1/2,\quad\beta=k_l-1/2.
\end{equation}
We identify $\lambda\in \a_\C^* $ with $(\lambda,2e_1)\in\C$ and 
let $t=e_1(\log a)/2$ be a coordinate on $A\simeq \R$. Then 
\begin{equation}
\rho(k)=k_{s}+2k_{l}=\alpha+\beta+1.
\end{equation}
The hypergeometric system (\ref{eqn:hgs}) turns out to be the differential 
equation 
\begin{equation}
\label{eqn:jacobi1}
L(k) F=(\lambda^2-\rho(k)^2)F, 
\end{equation}
where
\begin{equation}
\label{eqn:jde1}
L(k)=\frac{d^2}{dt^2}+2(k_{s}\coth t+2k_{l}\coth 2t)
\frac{d}{dt}
\end{equation}
and the hypergeometric function $F(\lambda,k;a_t)$ of type $BC_1$ 
is given by the Jacobi function
\begin{equation}
F(\lambda,k;a_t)=
\varphi^{(\alpha,\beta)}_{\sqrt{-1}\lambda}(t)
={}_2 F_1 \left(\tfrac12(\rho(k)-\lambda),\,
\tfrac12(\rho(k)+\lambda);\,\alpha+1;
-\sinh^2 t\right).
\end{equation}
Here ${}_2 F_1$ is the Gauss hypergeometric function. 
For $\lambda\not=1,2,\dots$, there is an 
another solution (\ref{eqn:hcseries}) 
of (\ref{eqn:jacobi1}) on $(0,\infty)$ given by
\begin{align}
\label{eqn:hcseries1}
\Phi&{}_{-\sqrt{-1}\lambda}^{(\alpha,\beta)}(t) \\
& =
(2\cosh t)^{\lambda-\rho(k)}{}_2F_1\left(
\tfrac12(\rho(k)-\lambda),\,\tfrac12(\alpha-\beta+1-\lambda);\,
1-\lambda;\,\cosh^{-2}t\right), \notag
\end{align}
which satisfies
\begin{equation}\label{eqn:asymhc}
\Phi_{-\sqrt{-1}\lambda}^{(\alpha,\beta)}(t)=e^{(\lambda-\rho)t}(1+o(t))
\text{ as }t\to\infty.
\end{equation}
For $\lambda\not\in\mathbb Z$ we have 
\begin{equation}
\varphi_{\sqrt{-1}\lambda}^{(\alpha,\beta)}(t)
=c_{\alpha,\beta}(-\sqrt{-1}\lambda)\Phi_{-\sqrt{-1}\lambda}^{(\alpha,\beta)}(t)
+
c_{\alpha,\beta}(\sqrt{-1}\lambda)\Phi_{\sqrt{-1}\lambda}^{(\alpha,\beta)}(t), 
\end{equation}
where
\begin{equation}
c_{\alpha,\beta}(-\sqrt{-1}\lambda)=c(\lambda,k)=
\frac{2^{\rho(k)-\lambda}\Gamma(\alpha+1)\Gamma(\lambda)}
{\Gamma(\frac12(\lambda+\rho(k)))\Gamma(\frac12(\lambda+\alpha-\beta+1))}.
\end{equation}

\section{Hypergeometric function of type $BC_n$}
\subsection{Commuting family of differential operators}
Let $n$ be a positive integer greater than $1$ and 
 $R$ be the root system of type $BC_n$
\begin{equation}
R_+=\{e_p,\,2e_p,\,e_i\pm e_j\,:\,1\leq p\leq n,\,1\leq i<j\leq n\}, 
\end{equation} 
where $\{e_1,\dots,e_n\}$ is the standard orthonormal basis 
of $E\simeq\R^n$. 
We call $\pm e_p,\,\pm(e_i\pm e_j)$, $\pm 2e_p$ short, middle, and 
long roots, respectively. 
We put 
\begin{equation}
k_{e_p}=k_s,\quad k_{e_i+e_j}=k_m,\quad k_{2e_p}=k_l
\end{equation}
for the multiplicities of short, middle, and long roots, respectively. 
Hereafter we assume that 
$k_m=0\text{ or }1$.
 Then the terms corresponding to 
the roots $e_i\pm e_j$ vanish in (\ref{eqn:laplace2}) and 
we have
\begin{align}
\label{eqn:laplace3}
\delta(k)^{\frac12} & 
\circ\{L(k)+(\rho(k),\rho(k))\}\circ\delta(k)^{-\frac12} \\
& = \sum_{j=1}^n 
\left(\partial_{e_j}^2+\frac{k_s (1-k_s-2k_l)}{(e^{\frac12e_j}-e^{-\frac12e_j})^2}+\frac{4k_l(1-k_l)}{(e^{e_j}-e^{-e_j})^2}\right).
\notag
\end{align}

Let $t_j=e_j(\log a)/2$ ($j=1,\dots,n$) be coordinates of $A\simeq \R^n$ and 
\begin{equation}
a_t=\exp(\textstyle\sum_{j=1}^n 2t_je_j).
\end{equation}
For $\lambda\in\h^*$ put $\lambda_j=(\lambda,2e_j)$. Then we have
\begin{equation}
\rho(k)_j=k_s+2k_l+2(n-j)k_m.
\end{equation}

Let $\Delta_m$ be the Weyl denominator associated with middle roots
\begin{align}
\Delta_m(a_t) &=
\prod_{\alpha\in R_+, \text{middle roots}}
(e^{\frac12\alpha}-e^{-\frac12\alpha}) \\
& =
2^{\frac12n(n-1)}
\prod_{1\leq i<j\leq n}(\cosh 2t_i-\cosh 2t_j).\notag
\end{align}
It is easy to see from (\ref{eqn:laplace3}) that 
\begin{equation}
\label{eqn:laplace4}
\Delta_m^{k_m}\circ (L(k)+(\rho(k),\rho(k))\circ \Delta_m^{-k_m}
=\sum_{j=1}^n L_j+n(k_s+2k_l)^2,
\end{equation}
where
\begin{equation}
\label{eqn:jde2}
L_j=\frac{\partial^2}{\partial t_j^2}+2
(k_{s}\coth t_j+2k_{l}\coth 2t_j)
\frac{\partial}{\partial t_j}.
\end{equation}

\begin{theorem}
\label{thm:main1}
If $k_{m}=0$ or $1$, then 
\[
\mathbb{D}(k)=\{D_p=\Delta_m^{-k_m}\circ p(L_1,\dots,L_n)\circ \Delta_m^{k_m}\,:\,
p \in \R[E]^W\}.
\]
 In particular, $\mathbb{D}(k)$ is 
generated by $D_{p_j}$ 
$(j=1,\dots,n)$,  
where $p_j$ is the $j$-th elementary symmetric function 
and $D_{p_1}=L(k)+(\rho(k),\rho(k))$.
\end{theorem}
\proof
Since $L_1,\dots,L_n$ mutually commute and 
\[
\gamma(k)(\Delta_m^{-k_m}\circ L_j\circ \Delta_m^{k_m})
=\partial_{e_j}^2-(k_s+2k_l)^2, 
\]
the theorem follows from Theorem~\ref{thm:ho1}.
\qed

\begin{remark}
{\em 
The right hand side of (\ref{eqn:laplace2}) has 
the form of a Schr\"odinger operator and 
Theorem~\ref{thm:ho1} tells that it defines a completely 
integrable system. 
Oshima \cite{Oshima} proved complete integrability of the 
Schr\"odinger operator
\[
P  =-\frac12\sum_{j=1}^n\partial_{e_j}^2+
\sum_{1\leq i<j\leq n}(u(t_i-t_j)+u(t_i+t_j))+\sum_{1\leq j\leq n}v(t_j)
\]
with
\begin{align*}
& u(x)=C_1\mathcal{P}(x)+C_2 \\
& v(x)=\frac{C_3 \mathcal{P}(x)^2+C_4 \mathcal{P}(x)^3+C_5 \mathcal{P}(x)^2
+C_6\mathcal{P}(x)+C_7}{\mathcal{P}'(x)^2}.
\end{align*}
If $C_1
=0$, then a result analogous to Theorem~\ref{thm:main1} 
holds. 
}
\end{remark}
\begin{remark}
{\em 
If $R$ is an arbitrary reduced root system and $k_\alpha=0$ or $1$ 
for 
all $\alpha\in R$, then the right hand side of (\ref{eqn:laplace2}) 
is just the Laplacian on the Euclidean space $E$. 
In this case, $\mathbb{D}(k)$ (taking conjugate by $\Delta_m^{k_m}$) consists 
of constant 
coefficient differential operators and the hypergeometric 
function is expressed by exponential functions. 
The case of all multiplicities equal to 1 
 is the case of 
 complex semisimple Lie groups in the sense of 
Remark~\ref{remark:sym}.  
Theorem~\ref{thm:main1} gives another case that $\mathbb{D}(k)$ has 
a simple simple expression. 
}
\end{remark}

\subsection{The hypergeometric function}

If $k_m=0$ or $1$, then the Harish-Chandra series 
(\ref{eqn:hcseries}) is given by a product of 
the Harish-Chandra series' of the form 
 (\ref{eqn:hcseries1}) for the root system of 
$R=BC_1$. 

\begin{proposition}\label{thm:hcseries}
Assume that $k_m=0$ or $1$ and 
let $\alpha=k_{s}+k_l-1/2,\,\beta=k_l-1/2$. 
If $\lambda$ satisfies condition 
{\rm (\ref{eqn:generic})}, then
\begin{equation}\label{eqn:hcseries-f}
\Phi(\lambda,k;a)=
{\Delta_m(a_t)}^{-k_m}\prod_{j=1}^n \Phi_{-\sqrt{-1}\lambda_j}^
{(\alpha,\beta)}(t_j).
\end{equation}
\end{proposition}
\proof
In view of (\ref{eqn:jde1}), (\ref{eqn:hcseries1}), 
(\ref{eqn:laplace4}), and (\ref{eqn:jde2}), 
the right hand side of (\ref{eqn:hcseries-f}) 
is a solution of (\ref{eqn:hgl}), where $\alpha$ and $\beta$ are 
given by (\ref{eqn:alpha})
.
We can see by elementary computations of power series 
that the right hand side of (\ref{eqn:hcseries-f})
 has a series expansion of the 
form (\ref{eqn:hcseries}) as in 
the same way as the proof of \cite[Theorem 1]{Hoogenboom}. 
By the uniqueness of the Harish-Chandra series, (\ref{eqn:hcseries}) 
follows. 
\qed\bigskip

By virtue of Proposition~\ref{thm:hcseries}, the hypergeometric 
function has a simple expression. 

\begin{theorem}\label{thm:main2}
Let $\alpha=k_{s}+k_l-1/2,\,\beta=k_l-1/2$  
and assume that $\alpha\not=0,-1,-2,\cdots$. 

If $k_m=1$, then 
\begin{equation}
\label{eqn:main1}
F(\lambda,k;a_t)=\frac{B}{\prod_{1\leq i<j\leq n}(\lambda_i^2-\lambda_j^2)}
\cdot
\frac{\det(\varphi^{(\alpha,\beta)}_{\sqrt{-1}\lambda_i}(t_j))_{1\leq i,\,j\leq n}}{\Delta_m (a_t)},
\end{equation}
where $B$ is given by 
\begin{equation}
\label{eqn:ca}
B=(-1)^{\frac12n(n-1)}2^{2n(n-1)}\prod_{i=1}^{n-1}((\alpha+i)^{n-i}i!).
\end{equation}
If $k_m=0$, then 
\begin{equation}
\label{eqn:main2}
F(\lambda,k;a_t)=
\frac{1}{n!}
\text{\rm perm}(\varphi^{(\alpha,\beta)}_{\sqrt{-1}\lambda_i}(t_j))_{1\leq i,\,j\leq n},
\end{equation}
where 
 $\text{\rm perm}(M)$ denotes the permanent 
$\sum_{\sigma\in S_n}m_1m_{\sigma(1)}\cdots m_n m_{\sigma(n)}$ 
of matrix $M=(m_{ij})_{1\leq i,j\leq n}$.
\end{theorem}
\proof
First notice that the Weyl group of type $BC_n$ is given by
\[
W=\{w=(\varepsilon,\sigma)\in \{-1\}^n\times S_n\,:\,
w(t_1,\dots,t_n)=(\varepsilon_1 t_{\sigma(1)},\dots,
\varepsilon_n t_{\sigma(n)})\}.
\]

Assume that $k_m=1$.  
The $c$-function for the middle roots (the product 
is take over the middle roots in (\ref{eqn:cf1})) is given 
by 
\begin{align}
\tilde{c}_m (\lambda,k) & =\prod_{1\leq i<j\leq n}
\frac{\Gamma\left(\frac12(\lambda_i+\lambda_j)\right)
\Gamma\left(\frac12(\lambda_i-\lambda_j)\right)}
{\Gamma\left(\frac12(\lambda_i+\lambda_j)+1\right)
\Gamma\left(\frac12(\lambda_i-\lambda_j)+1\right)} \\
& =
\frac{2^{n(n-1)}}{\prod_{1\leq i<j\leq n}(\lambda_i^2-\lambda_j^2)}.
\notag
\end{align}
The $c$-function for $e_j$ and $2e_j$ is given by
\begin{align}
\tilde{c}_{e_j}(\lambda,k)\tilde{c}_{2e_j}(\lambda,k)
& = 
\frac{2^{-\lambda_j-k_s+1}\Gamma(\lambda_j)}
{\Gamma\left(\frac12(\lambda_j+k_s+1)\right)
\Gamma\left(\frac12(\lambda_j+k_s+2k_l)\right)} \\
& = 
2^{-2k_s-2k_l+1}\Gamma\left(k_s+k_l+\tfrac12\right)^{-1}
c_{\alpha,\beta}(-\sqrt{-1}\lambda_j).\notag
\end{align}
We have
\begin{align}
\tilde{c}(\lambda,k)
& =\tilde{c}_m(\lambda,k)\prod_{j=1}^n 
\tilde{c}_{e_j}(\lambda,k)\tilde{c}_{2e_j}(\lambda,k) \\
& =
\frac{2^{n(n-2k_s-2k_l)}}
{{\Gamma\left(k_s+k_l+\tfrac12\right)^n}
\prod_{1\leq i<j\leq n}(\lambda_i^2-\lambda_j^2)}
\prod_{j=1}^n 
c_{\alpha,\beta}(\lambda_j).
\notag
\end{align}

The hypergeometric function is given by
\begin{align*}
\Delta_m & (a_t)F(\lambda,k;a_t)  = 
\tilde{c}(\rho(k),k)^{-1}
\sum_{w\in W}\tilde{c}(w\lambda,k)
\Delta_m(a_t)
\Phi(w\lambda,k,a_t)
 \\
& = 
B
\sum_{\sigma\in S_n,\,\varepsilon\in\{-1\}^n}
\frac{1}{\prod_{i<j}(\lambda_{\sigma(i)}^2-\lambda_{\sigma(j)}^2)}
\prod_{l=1}^n
c_{\alpha,\beta}(-\sqrt{-1}\varepsilon_l\lambda_{\sigma(l)})
\Phi^{(\alpha,\beta)}_{-\sqrt{-1}\varepsilon_l\lambda_{\sigma(l)}}(t_l) \\
& = 
B
\frac{1}{\prod_{i<j}(\lambda_{i}^2-\lambda_{j}^2)}
\sum_{\sigma\in S_n}{\text{sgn}\,\sigma}
\prod_{l=1}^n
\varphi_{\sqrt{-1}\lambda_{\sigma(l)}}^{(\alpha,\beta)}(t_l) \\
&= B
\frac{\det(\varphi_{\sqrt{-1}\lambda_{i}}^{(\alpha,\beta)}(t_j))_{i,j}}
{\prod_{i<j}(\lambda_{i}^2-\lambda_{j}^2)},
\end{align*}
where
\[
B=\frac{2^{n(n-1)}}{\tilde{c}(\rho(k),k)
\left(2^{2k_s+2k_l-1}\Gamma\left(k_s+k_l+\tfrac12\right)\right)^n}.
\]
The formula for $B$ can be obtained by explicit computations. 

Next suppose $k_m=0$. Then $c_m(\lambda,k)=
\lim_{k_m\to 0}\tilde{c}(\lambda,k)/\tilde{c}(\rho(k),k)=1/n!$. 
Here $c_m(\lambda,k)$ is the $c$-function for the middle 
roots (the product 
is take over the middle roots in (\ref{eqn:cf2})). 
(\ref{eqn:main2}) follows by direct computation 
similar to that of deriving (\ref{eqn:main1}). 

\qed

\begin{remark}\label{rem:group1}
{
\rm 
Let $p$ and $q$ ($p\leqq q$) be positive integers and
put $k_s=q-p$, $k_m=1$, and $k_l=1/2$. 
Then the hypergeometric function $F(\lambda,k;a_t)$ is 
the radial part of the spherical function on 
$SU(p,q)/S(U(p)\times U(q))$. 
In this case Theorem~\ref{thm:main1}, Theorem~\ref{thm:hcseries}, 
and Theorem~\ref{thm:main1} were given by 
Berezin and Karpelevi\v c \cite{BK} without proof 
and a complete proof was 
given by Hoogenboom~\cite{Hoogenboom}. 
}
\end{remark}

We give two corollary of our results. 

First we give a limit case of the hypergeometric function. 
We replace $(t,\lambda)$ by $(\epsilon t,\epsilon^{-1}\lambda)$ 
and let $\epsilon\downarrow 0$. Then the hypergeometric 
equation (\ref{eqn:jacobi1}) of type $BC_1$ becomes 
\begin{equation}
\label{eqn:bessel1}
\frac{d^2 u}{dt^2}+\frac{2\alpha+1}{t}\frac{\partial ^2 u}{dt^2}=\lambda^2 u.
\end{equation}
Here we put $\alpha=k_s+k_l+1/2$. 
There exists a unique even solution  of (\ref{eqn:bessel1}) that is 
regular at $0$ and $u(0)=1$, which is given by
\begin{equation}
\label{eqn:bessel2}
\mathcal{J}_\alpha(\sqrt{-1}\lambda t)=
2^\alpha\Gamma(\alpha+1)(\sqrt{-1}\lambda t)^{-\alpha}J_\alpha(\sqrt{-1}\lambda t),
\end{equation}
where $J_\alpha$ denote the usual Bessel function. Then 
it is known  \cite[\S 2.3]{Koornwinder} that 
\begin{equation}
\label{eqn:limrank1}
\lim_{\epsilon\downarrow 0}\varphi^{(\alpha,\beta)}_{\sqrt{-1}\epsilon^{-1}\lambda}
(\epsilon t)=\mathcal{J}_\alpha(\sqrt{-1}\lambda t).
\end{equation}

The limit of operator (\ref{eqn:laplace1}) become
\begin{equation}
L(k)^\text{rat}=\sum_{j=1}^n \partial_{\xi_j}^2+\sum_{\alpha\in R_+}
\frac{2k_\alpha}{\alpha}
\partial_\alpha
\end{equation}
and we have 
\begin{equation}
\lim_{\epsilon\downarrow 0}\epsilon^{-n(n-1)}\Delta_m(a_{\epsilon t})
=\prod_{\alpha\in R_+,\text{middle roots}}\alpha(\log a_t). 
\end{equation}
We denote the right hand side of the above equation by $\Delta_{m,\text{rat}}(a_t)$. 
Put 
\begin{equation}
L^\text{rat}_j=\frac{\partial^2}{\partial t_j^2}+\frac{2k_s+2k_l+2}{t_j}\frac{\partial}{\partial t_j}. 
\end{equation}
Then we have the following explicit expression of commuting family of 
differential operators including $L(k)^\text{rat}$. 

\begin{corollary}
\label{cor:main1}
If $k_{m}=0$ or $1$, then 
\[
\{D_p^\text{\rm rat}=
\Delta_{m,\text{\rm rat}}^{-k_m}\circ 
p(L_1,\dots,L_n)\circ \Delta_{m,\text{\rm rat}}^{k_m}\,:\,
p \in \R[E]^W\}
\]
forms a commutative algebra of differential operators, which is 
generated by $\Delta_{m,\text{\rm rat}}^{-k_m}\circ p_j(L_1,\dots,L_n)\circ
 \Delta_{m,\text{\rm rat}}^{k_m}$, 
$(j=1,\dots,n)$, 
where $p_j$ is the $j$-th elementary symmetric function. 
$D_{p_1}^\text{\rm rat}=L(k)^{\text{\rm rat}}$ and 
the principal symbol of $D_{p_j}^\text{\rm rat}$ is $p_j$ for $j=1,\dots,n$. 
\end{corollary}

By Theorem~\ref{thm:main2} and (\ref{eqn:limrank1}) we have the 
following limit formula.

\begin{corollary}
\label{cor:main2}
Let $\alpha=k_{s}+k_l-1/2$   
and assume that $\alpha\not=0,-1,-2,\cdots$ and $\lambda_j\not=0$,  
$t_j\not=0$ $(j=1,\dots,n)$. 

If $k_m=1$, then 
\begin{equation}
\label{eqn:main3}
\lim_{\epsilon\downarrow 0}
F(\epsilon^{-1}\lambda,k;a_{\epsilon t})=
\frac{A}{\prod_{1\leq i<j\leq n}(\lambda_i^2-\lambda_j^2)}
\cdot
\frac{\det(\mathcal{J}_\alpha({\sqrt{-1}\lambda_i}t_j))_{1\leq i,\,j\leq n}}
{\Delta_{m,\text{\rm rat} }(a_t)},
\end{equation}
where $A$ is given by (\ref{eqn:ca}). 
If $k_m=0$, then 
\begin{equation}
\label{eqn:main4}
\lim_{\epsilon\downarrow 0}
F(\epsilon^{-1}\lambda,k;a_{\epsilon t})=
\frac{1}{n!}
\text{\rm perm}(\mathcal{J}_\alpha(\sqrt{-1}\lambda_i t_j))_{1\leq i,\,j\leq n}.
\end{equation}
\end{corollary}

\begin{remark}
{\rm 
In the group case that we mentioned in Remark~\ref{rem:group1}, 
(\ref{eqn:main3}) was proved by Meaney~\cite{M}. 
It gives contraction of spherical functions between 
symmetric spaces of the non-compact type and the Euclidean type. 

The right hand side of (\ref{eqn:main3}) and (\ref{eqn:main4}) give 
explicit expression for   
the Bessel function of type $BC_n$ which was defined by 
Opdam~\cite[Definition 6.9]{Opdam}. The Bessel function of type $BC_n$ 
for $k_m=0$ or $1$ is a $W$-invariant $C^\infty$ joint-eigenfunction of 
the commuting family of differential operators given in 
Corollary~\ref{cor:main1} being equal to 1 at the origin. 

The type of limit transition in Corollary~\ref{cor:main2} was given 
also by Ben Sa\"{i}d and \O rsted \cite{BO1,BO2}, and 
de Jeu~\cite{deJeu}. 
}
\end{remark}

Finally we give a formula for a $\Theta$-spherical function. 
Let $\Psi$ denote the set of simple roots in $R_+$, 
\[
\Psi=\{e_1-e_2,\dots,e_{n-1}-e_n,e_n\}.
\]
For a subset $\Theta\subset \Psi$, let 
$\langle \Theta\rangle=R\cap \sum_{\alpha\in\Theta} 
\mathbb Z\alpha$ and define 
$\tilde{c}_\Theta(\lambda,k)$ by the product 
of the form 
(\ref{eqn:cf1}) where the product is taken over 
$R_+\cap\langle\Theta\rangle$ and let 
$c_\Theta(\lambda,k)=\tilde{c}_\Theta(\lambda,k)/
\tilde{c}_\Theta(\rho(k),k)$. 

We make a sum 
\begin{equation}
\label{eqn:partialhg}
F_\Theta(\lambda,k,a)=\sum_{w\in W_\Theta}
c_\Theta(w\lambda,k)\Phi(w\lambda,k;a). 
\end{equation}
The sum of the form 
(\ref{eqn:partialhg}) is important in harmonic analysis 
of the spherical function on symmetric spaces 
(c.f. \cite{OP}, \cite[Chapter 6]{Sc}, 
\cite{Sh}). 

By Theorem~\ref{thm:hcseries}, we can derive 
formulae for $F_\Theta(\lambda,k,a)$. 
For $\Theta=\Psi\setminus\{e_{1}-e_2,\dots,e_{j-1}-e_j\}\,
(2\leq j\leq n)$ we have a formula for $F_\Theta(\lambda,k;a_t)$ 
that is similar to the formula for $F(\lambda,k;a_t)$ in 
 Theorem~\ref{thm:main2}. 

If 
 $\Theta=\{e_1-e_2,\dots,e_{n-1}-e_n\}$, 
then $\langle\Theta\rangle$ is a root system of type $A_{n-1}$ and 
we have the following result. 

\begin{corollary}\label{thm:hcseries-2}
Assume that $k_m=0$ or $1$ and 
let $\Theta=\{e_1-e_2,\dots,e_{n-1}-e_n\}$ and 
$\alpha=k_{s}+k_l-1/2,\,\beta=k_l-1/2$. 
Then $F_\Theta(\lambda,k;a_t)$ is holomorphic in $\lambda$ in 
the region $\text{\rm Re}\, 
\lambda_i>0\,(i=1,\dots,n)$. 
Moreover we have the following results.
\\
{\rm (i)} Suppose $k_m=1$ and put 
 $\pi(x_1,\dots,x_n)=\prod_{1\leq i<j\leq n}(x_i-x_j)$.  
Then we have
\begin{equation}
\label{eqn:main2-1}
F_\Theta(\lambda,k;a_t)=
\frac{\pi(
\rho
(k)
)
}{\pi(\lambda)}
\cdot
\frac{\det(\Phi^{(\alpha,\beta)}_{\sqrt{-1}\lambda_i}
(t_j))_{1\leq i,\,j\leq n}}{\Delta_m (a_t)}.
\end{equation}
Moreover, if $\text{\rm Re}\, 
\lambda_i>0\,(i=1,\dots,n)$, then 
\begin{equation}
\label{eqn:main2-2}
\lim_{u\to \infty}e^{(\rho(k)-\lambda)(\log \,a_{(u,\dots,u)})}
F_\Theta(\lambda,k;a_{(t_1+u,\dots,t_n+u)})=
\frac{\pi(\rho(k))}{\pi(\lambda)}
\cdot
\frac{\det(e^{\lambda_i t_j})_{1\leq i,\,j\leq n}}
{\pi(e^{2t_1},\dots,e^{2t_n})}.
\end{equation}
\\
{\rm (ii)} 
If $k_m=0$, then 
\begin{equation}
\label{eqn:main2-3}
F_\Theta(\lambda,k;a_t)=\frac{1}{n!}
\text{\rm perm}(\Phi^{(\alpha,\beta)}_{\sqrt{-1}\lambda_i}
(t_j))_{1\leq i,\,j\leq n}.
\end{equation}
Moreover, if $\text{\rm Re}\, 
\lambda_i>0\,(i=1,\dots,n)$, then 
\begin{equation}
\label{eqn:main2-4}
\lim_{u\to \infty}e^{(\rho(k)-\lambda)(\log \,a_{(u,\dots,u)})}
F_\Theta(\lambda,k;a_{(t_1+u,\dots,t_n+u)})=
\frac{1}{n!}
{\text{\rm perm}(e^{\lambda_i t_j})_{1\leq i,\,j\leq n}}.
\end{equation}
\end{corollary}

\proof
$F_\Theta(\lambda,k;a_t)$ is holomorphic in 
the region $\text{\rm Re}\, 
\lambda_i>0\,(i=1,\dots,n)$ by \cite[Theorem 8]{OP}. 
(\ref{eqn:main2-1}) and (\ref{eqn:main2-3}) follows by 
simple computations.
(\ref{eqn:main2-2}) and (\ref{eqn:main2-4}) follows from 
(\ref{eqn:asymhc}). 
\qed

\begin{remark}{\rm (i) 
The right hand sides of 
(\ref{eqn:main2-2}) and (\ref{eqn:main2-4}) are  
hypergeometric function of type $A_{n-1}$ with 
the multiplicity $1$ and $0$ respectively. 
Namely, the right hand side of (\ref{eqn:main2-2}) is 
the spherical function on $SL(n,\mathbb C)/SU(n)$ 
(c.f. \cite[Chapter IV Theorem 5.7]{Helgason}) and 
(\ref{eqn:main2-4}) is the normalized average of the 
exponential function $e^{(\lambda,t)}$ under the 
action of the symmetric group. 

\smallskip
\noindent
(ii) By \cite[Proposition~2.6, Remark~6.13]{Sh}, the 
spherical function for a one-dimensional $K$-type $(\tau_{-\ell_1},\tau_{-\ell_2})$ 
on $SU(p,q)$ can be written as the 
hypergeometric function $F(\lambda,k;a_t)$ with 
$k_s=m/2-\ell_2,\,k_m=1,\,k_l=1/2-\ell_1-\ell_2$. 
Here $m=1$ and $\ell_1=\ell_2$  if $p\not=q$, and 
$m=0$ if $p=q$. Thus spherical functions for one-dimensional $K$ types 
on  $SU(p,q)$ are given by Theorem~\ref{thm:main1}. 
Conversely, by considering the universal covering group of  $SU(p,q)$, 
we can take $\ell_1,\,\ell_2$ arbitrary complex numbers, hence the 
hypergeometric function (\ref{eqn:main1}) for any $k_s$ and $k_l$ 
corresponds to a spherical function on $\widetilde{SU(p,p)}$. 

By the above observation, 
the Plancherel formula for the integral transform with the kernel 
$F(\lambda,k;a)$ with $k_m=1$ 
is a special case of \cite[Theorem~6.11]{Sh}. Notice that 
low dimensional spectra including discrete spectra appear in general. 
It seems to be possible to give an alternative proof of the Plancherel formula by 
rank one reduction as in \cite[Theorem 22]{M}. 

\smallskip
\noindent
(iii) In Theorem~\ref{thm:main1} we give an explicit formula for the hypergeometric function 
of type $BC_n$ with $k_m=0,\,1$ and $k_s,\,k_l$ arbitrary. 
We obtain a formula of the hypergeometric function for $k_m\in\mathbb{Z}$ by applying 
Opdam's 
hypergeometric shift operator corresponding to the middle roots, 
which is a differential operator of order $n(n-1)/2$ (c.f. \cite[Definition 3.2.1]{Heckman}). 
}
\end{remark}

\end{document}